\documentclass[seceqn,dvips]{arxbj}
\usepackage{graphics}
\aid{0}
\volume{16}
\issue{2}
\pubyear{2010}
\firstpage{493}
\lastpage{513}
\doi{10.3150/09-BEJ219}
\makeatletter
\def\epsilon{\varepsilon}

\renewcommand{\citep}[1]{(\citet{#1})}

\newtheorem{proposition}{Proposition}
\newtheorem{theorem}{Theorem}
\newtheorem{corollary}{Corollary}
\newtheorem{lemma}{Lemma}

\makeatother

\begin{document}
\begin{frontmatter}

\title{Estimation of a probability with optimum guaranteed confidence
in inverse binomial sampling}
\runtitle{Optimum guaranteed confidence in inverse binomial sampling}

\begin{aug}
\author{\fnms{Luis} \snm{Mendo}\thanksref{e1}\ead[label=e1,mark]{lmendo@grc.ssr.upm.es}\corref{}}
\and
\author{\fnms{Jos\'{e} M.} \snm{Hernando}\thanksref{e2}\ead[label=e2,mark]{hernando@grc.ssr.upm.es}}
\pdfauthor{Luis Mendo, Jose M. Hernando}
\runauthor{L. Mendo and J.M. Hernando}
\address{E.T.S. Ingenieros de Telecomunicaci\'{o}n, Polytechnic
University of Madrid,
28040 Madrid, Spain.\\
E-mails: \printead*{e1,e2}}
\end{aug}

% HISTORY:
\received{\smonth{9} \syear{2008}}

% ABSTRACT
%
\begin{abstract}
Sequential estimation of a probability $p$ by means of inverse binomial
sampling is considered. For $\mu_1,\mu_2>1$ given, the accuracy of an
estimator $\hat p$ is measured by the confidence level $P[p/\mu_2
\leq\hat p \leq p \mu_1]$. The confidence levels $c_0$ that can be
guaranteed for $p$ unknown, that is, such that $P[p/\mu_2 \leq
\hat p \leq p\mu_1] \geq c_0$ for all $p \in(0,1)$, are investigated.
It is shown that within the general class of randomized or
non-randomized estimators based on inverse binomial sampling, there is
a maximum $c_0$ that can be guaranteed for arbitrary~$p$. A
non-randomized estimator is given that achieves this maximum guaranteed
confidence under mild conditions on $\mu_1$, $\mu_2$.
\end{abstract}

% KEYWORDS
%

\end{frontmatter}

%s1 ###
\section{Introduction}
\label{parte:intro}

In a sequence of Bernoulli trials with probability of success $p$ at
each trial, consider the estimation of $p$ by inverse binomial
sampling. This sampling scheme, first discussed by \citet{Haldane45},
consists in observing the sequence until a given number $r$ of
successes is obtained. The resulting number of trials $N$ is a
sufficient statistic for $p$ (\citet{Lehmann98}, page~101). The uniformly
minimum variance unbiased estimator of $p$ is (\citet{Mikulski76})
%
%e1.1 ###
\begin{equation}\label{eq:hatpN-1n-1}
\hat p= \frac{r-1}{N-1}
\end{equation}
and for $r\geq3$, it has a normalized mean square error $\mathrm{E}
[(\hat p-p)^2]/p^2$ (or $\operatorname{Var}[\hat p]/p^2$) lower than $1/(r
-2)$, irrespective of $p$ (\citet{Mikulski76}; \citet{Sathe77}; \citet{Prasad82}).

This paper analyzes inverse binomial sampling from a different point of
view, related to interval estimation. Given $\mu_1, \mu_2> 1$, the
accuracy of an estimator $\hat p$ is measured by the probability~$c$
that $\hat p$ lies in the interval $[p/\mu_2, p\mu_1]$. The
motivation to use a relative interval $[p/\mu_2, p\mu_1]$, instead of
an interval $[x_1,x_2]$ with $x_1, x_2$ fixed, is the fact that $P
[p/\mu_2\leq\hat p\leq p\mu_1]$, unlike $P[x_1 \leq\hat p
\leq x_2]$, has a definite meaning independent of $p$. Moreover, using
a relative interval allows another interpretation of $c$: since
$p/\mu_2\leq\hat p\leq p\mu_1$ if and only if $\hat p/\mu_1\leq
p \leq\hat p\mu_2$, $c$ gives the probability that the true value
$p$ is covered by the random interval $[\hat p/\mu_1, \hat p\mu_2
]$. In the sequel, $c$ will be referred to as the \textit{confidence}
(or confidence level) associated with $\mu_1$, $\mu_2$.

Recently, \citeauthor{Mendo06} (\citeyear{Mendo06,Mendo08a}) have shown that the estimator
%
%e1.2 ###
\begin{equation}
\label{eq:hatpN-1n}
\hat p= \frac{r-1}{N}
\end{equation}
has the following properties for $r\geq3$. The confidence level
$c$ associated with $\mu_1$, $\mu_2$ has an asymptotic value $\bar c$
as $p \rightarrow0$, namely
$
\bar c = \gamma(r,(r-1)\mu_2)- \gamma(r,(r
-1)/\mu_1),
$
where
\[
\gamma(r,t) = \frac{1}{\Gamma(r)} \int_0^t s^{r
-1}\exp(-s)\,\mathrm{d}s
\]
is the regularized incomplete gamma function. Furthermore, the
confidence for any $p \in(0,1)$ exceeds this asymptotic value provided
that $\mu_1$ and $\mu_2$ satisfy certain lower bounds.
Similar results have been established \citep{Mendo08b} for the
uniformly minimum variance unbiased estimator (\ref{eq:hatpN-1n-1}).
Using a more general setting in which the Bernoulli random variables
are replaced by arbitrary bounded random variables, \citet{Chen07env}
has obtained comparable (although somewhat less tight) results for
(\ref{eq:hatpN-1n-1}), as well as for the maximum likelihood estimator
%
%e1.3 ###
\begin{equation}
\label{eq:hatpNn}
\hat p= \frac{r}{N}.
\end{equation}

From the aforementioned results, the question naturally arises as to
whether the attained confidence could be improved using estimators
other than
(\ref{eq:hatpN-1n-1})--(\ref{eq:hatpNn}).
% *!* --
This motivates the study of arbitrary estimators based on inverse
binomial sampling. In this regard, it is noted that $1-c$ can be
expressed as the risk $\mathrm{E}[L(\hat p)]$ corresponding to
the loss
function $L$ defined as
%
%e1.4 ###
\begin{equation}
\label{eq:lossfunction}
L(x) =
\cases{
0, &\quad if $ x \in[p/\mu_2, p\mu_1]$, \cr
1, &\quad otherwise.
}
\end{equation}
Since $N$ is a sufficient statistic, for any estimator defined
in terms of the observed Bernoulli random variables, there exists an
estimator that depends on the observations through $N$ only and
that has the same risk; this equivalent estimator is possibly a
randomized one (\citet{Lehmann98}, page 33). (The fact that $L$ is
non-convex prevents application of a corollary of the Rao--Blackwell
theorem (\citet{Lehmann98}, page 48) to discard randomized estimators.)
Thus, attention can be restricted to estimators that depend on the
observed variables only through~$N$; however, both randomized
and non-randomized estimators have to be considered.

The purpose of the present paper is to investigate the limits on the
confidence levels $c_0$ that can be guaranteed (in the sense that the
actual confidence equals or exceeds $c_0$, irrespective of~$p$) in
inverse binomial sampling. More specifically, the objectives are:
\begin{itemize}
\item
to determine the supremum of $\inf_{p \in(0,1)} c$ over all
(randomized or non-randomized) estimators based on inverse binomial sampling;
\item
to find an estimator, if it exists, that can achieve this supremum.
\end{itemize}
Consequently, the main focus of this paper is on non-asymptotic
results, valid for $p \in(0,1)$. Nonetheless, asymptotic results for
$p \rightarrow0$ will also be derived since, apart from their own
theoretical importance, they provide an upper bound on (and an
indication of) what can be achieved for arbitrary $p$.

Section~\ref{parte:asympt} shows that $\liminf_{p\rightarrow0}c$
has a maximum over all estimators based on inverse binomial sampling
and computes this maximum. This sets an upper bound on the confidence
that can be guaranteed by an arbitrary estimator. Section~\ref{parte:noasympt} establishes that this upper bound is also a maximum, that
is, estimators exist that can achieve this guaranteed confidence.
Specifically, estimators are given that guarantee the maximum
confidence for sufficiently small $p$ and that guarantee the maximum
confidence for any $p$, under certain conditions on the relative
interval being considered. Section \ref{parte:disc} discusses the
results, comparing them with those from other works. Section \ref
{parte:proofs} contains proofs of all results in the paper.

%s2 ###
\section{Asymptotic analysis}
\label{parte:asympt}

It is assumed in the sequel that $r\geq3$. Let $t^{(j)}$
denote $t(t-1)\cdots(t-j+1)$. The probability function of $N$,
$f(n) = P[N= n]$, is
%
%e2.1 ###
\begin{equation}
\label{eq:f}
f(n)
= \frac{(n-1)^{(r-1)}}{(r-1)!} p^r(1-p)^{n
-r} \qquad\mbox{for }n\geq r.
\end{equation}

As justified in Section \ref{parte:intro}, it is sufficient to
consider (possibly randomized) estimators defined in terms of the
sufficient statistic $N$. Let $\mathcal{F}$ denote the set of all
functions from $\{r, r+1, r+2,\ldots\}$ to $\mathbb R$.
A \textit{non-randomized estimator} $\hat p$ can be described as
$\hat p= g(N)$, with $g\in\mathcal{F}$. Thus, a
non-randomized estimator is entirely specified by its function $g
$. A \textit{randomized estimator} is a random variable $\hat p$
whose distribution depends, in general, on the value taken by $N
$. Let $\Pi_n$ denote the distribution function of $\hat p$
conditioned on $N= n$.
% that is, $\Pi_\nden(x) = \Pr[\hatvap\leq x \st\vanden= \nden]$.
The randomized estimator is completely specified by the functions $\Pi
_n$, $n\geq r$. Thus, denoting by $\mathcal{F}_{\mathcal{R}}$ the
class of
all functions from $\{r, r+1, r+2,\ldots\}$ to the set
of real functions of a real variable, a randomized estimator is defined
by a function $G\in\mathcal{F}_{\mathcal{R}}$ that assigns
%the distribution function
$\Pi_n$ to each $n$.

Non-randomized estimators form a subset of the class of randomized estimators.
% ; namely, for $\fest\in\Fe$, the corresponding $\festa\in\FeA$ is
%defined by the property that $\Pi_\nden(x)=1$ if $x \geq\fest(\nden)$
%and $0$ otherwise.
Thus, any statement that applies to randomized estimators will also be
valid, in particular, for non-randomized estimators. The reason that
the specialized class of non-randomized estimators has been explicitly
defined is that, on one hand, their simplicity makes them more
attractive for applications and, on the other hand, it will be seen
that, under certain conditions, no loss of optimality in guaranteed
confidence is incurred by restricting to non-randomized estimators.
Throughout the paper, when referring to an arbitrary estimator without
specifying its type, the general class of randomized estimators
(including the non-randomized ones) will be meant.

Given any estimator $G\in\mathcal{F}_{\mathcal{R}}$, the confidence
associated with
$\mu_1$, $\mu_2$ will be expressed in the sequel as a function $c(p)$
(that is, the dependence on $p$ will be explicitly indicated). Given
$r$, $\mu_1$ and $\mu_2$, the latter function is determined by
$G$. An estimator is said to \textit{guarantee} a confidence
level $c_0$ in the interval $(p_1,p_2)$ if $c(p) \geq c_0$ for all $p
\in(p_1,p_2)$. If the interval is $(0,1)$, then the estimator is said
to \textit{globally guarantee} this confidence level. An estimator
\textit{asymptotically guarantees} a certain $c_0$ if there exists
$\epsilon> 0$ such that the estimator guarantees $c_0$ in the interval
$(0,\epsilon)$.
%Global guarantee implies asymptotic guarantee.

The problem being addressed can be rephrased as that of optimizing the
globally guaranteed confidence within the general class of estimators
based on inverse binomial sampling, or finding a minimax estimator with
respect to the risk defined by the loss function (\ref{eq:lossfunction}). The following proposition, and its ensuing
particularization, provides the motivation for studying $\liminf_{p
\rightarrow0} c(p)$ as part of this optimization problem.

\begin{proposition}
\label{prop:asintnoasint}
If a given estimator, with confidence function $c(p)$, guarantees a
confidence level $c_0$ in an interval $(p_1,p_2)$, then, necessarily,
$\liminf_{p \rightarrow p_0} c(p) \geq c_0$ for any $p_0 \in[p_1, p_2]$.
\end{proposition}

Particularizing Proposition \ref{prop:asintnoasint} to
% vale
$p_1=p_0=0$, it is seen that $\liminf_{p \rightarrow0} c(p)$
represents an upper bound on the confidence levels that can be globally
or asymptotically guaranteed.

An important subclass of non-randomized estimators is formed by those
defined by functions $g\in\mathcal{F}$ for which $c(p)$ has an
asymptotic value, that is, for which $\lim_{p \rightarrow0} P
[p/\mu_2\leq g(N) \leq p\mu_1]$ exists. The set of all
such functions will be denoted $\mathcal{F}_{\mathrm{p}}$. As will be
seen, another
important subclass is that corresponding to the set of functions $g
\in\mathcal{F}$ for which $\lim_{n\rightarrow\infty} ng
(n)$ exists, is finite and non-zero. This set will be denoted
$\mathcal{F}_{\mathrm{n}}$. The following result establishes that
$\mathcal{F}_{\mathrm{n}}\subset\mathcal{F}_{\mathrm{p}}$.

\begin{proposition}
\label{prop:limdeconf}
For a non-randomized estimator defined by $g\in\mathcal{F}_{\mathrm
{n}}$, $\lim_{p
\rightarrow0} c(p)$ exists and equals $\bar c$, given by
%
%e2.2 ###
\begin{equation}
\label{eq:barc}
\bar c = \gamma(r,\Omega\mu_2) - \gamma(r,\Omega/\mu_1
),\qquad\Omega= \lim_{n\rightarrow\infty} ng(n).
\end{equation}
\end{proposition}

The converse of Proposition \ref{prop:limdeconf} is not true; that
is, $\mathcal{F}_{\mathrm{n}}\neq\mathcal{F}_{\mathrm{p}}$.
% a function $\fest\in\Felp$ does not necessarily belong to $\Feln$.
A simple counterexample is given by
\[
g(n) =
\cases{
{\omega'}/{n}, & \quad if there exists $k \in\mathbb N $
such that $ n= 2^k$, \cr
{\omega}/{n}, &\quad otherwise,
}
\]
with $\omega' \neq\omega$. It is easily seen that this function is
in $\mathcal{F}_{\mathrm{p}}$, with $\lim_{p \rightarrow0} c(p) =
\gamma(r,\omega
\mu_2) - \gamma(r,\omega/\mu_1)$; however, $\lim_{n
\rightarrow\infty} ng(n)$ does not exist.

Given $r$, $\mu_1$ and $\mu_2$, the maximum of $\bar c$ over all
$g\in\mathcal{F}_{\mathrm{n}}$ is attained when $\lim_{n\rightarrow
\infty}
ng(n)$ equals $\Omega^*$, given by
%
%e2.3 ###
\begin{equation}\label{eq:Omega*}
\Omega^* = r  \frac{\log\mu_2- \log(1/\mu_1)}{\mu_2
-1/\mu_1},
\end{equation}
as is readily seen by differentiating $\bar c$ in (\ref{eq:barc}).
Let $c^*$ denote the resulting maximum.
%c^* = \gamma(\nnum,\Omega^* \fdos) - \gamma(\nnum,\Omega^*/\funo).
Defining
%
%e2.4 ###
\begin{equation}
\label{eq:fud}
M= \mu_1\mu_2,
\end{equation}
the terms $\Omega^*/\mu_1$ and $\Omega^*\mu_2$ can be expressed as
%
%e2.5 ###
\begin{equation}
\label{eq:Omega*funofdos}
\Omega^*/\mu_1= \frac{r\log M}{M-1}, \qquad
\Omega^*\mu_2= \frac{rM\log M}{M-1}
\end{equation}
and thus
%
%e2.6 ###
\begin{equation}
\label{eq:c*fud}
c^* = \gamma \biggl(r, \frac{rM\log M}{M-1}
 \biggr)
- \gamma \biggl(r, \frac{r\log M}{M-1}  \biggr).
\end{equation}
The following theorem establishes that $c^*$ is not only the maximum of
$\lim_{p \rightarrow0} c(p)$ within the subclass of non-randomized
estimators defined by $\mathcal{F}_{\mathrm{n}}$, but also the
maximum of $\liminf_{p
\rightarrow0} c(p)$ within the general class of randomized estimators
defined by $\mathcal{F}_{\mathcal{R}}$.
%Thus, optimality of $\liminf_{p \rightarrow0} c(p)$ is not lost by
%restricting to non-randomized estimators.

\begin{theorem}\label{theo:caracdemaxasint}
The maximum of $\liminf_{p \rightarrow0} c(p)$ over all estimators
defined by functions $G\in\mathcal{F}_{\mathcal{R}}$ is $c^*$, given
by (\ref{eq:c*fud}).
\end{theorem}

As can be seen from (\ref{eq:c*fud}), $c^*$ depends on $\mu_1$ and
$\mu_2$ only through $M$. An explanation of this result is as
follows. Given $\mu_1$, $\mu_2$, let $\hat p$ be an arbitrary
estimator and consider $a>0$. If $\mu_1$ and $\mu_2$ are replaced by
$\mu_1'=a\mu_1$ and $\mu_2'=\mu_2/a$, respectively, defining a
modified estimator $\hat p' = a \hat p$, it is clear that $p/\mu_2
' \leq\hat p' \leq p\mu_1'$ if and only if $p/\mu_2\leq\hat p
\leq p\mu_1$. This shows that any value of $\liminf_{p\rightarrow0}
c(p)$ that can be achieved for $\mu_1$, $\mu_2$ can also be achieved
for $a\mu_1$, $\mu_2/a$ (using a different estimator), and
conversely. Thus, $c^*$ is the same for $\mu_1$, $\mu_2$ and for
$a\mu_1$, $\mu_2/a$.

%s3 ###
\section{An optimum estimator for certain relative intervals}
\label{parte:noasympt}

The asymptotic results in Section \ref{parte:asympt} impose a limit
on the confidence levels that can be guaranteed, as established by the
following corollary of Proposition \ref{prop:asintnoasint} and
Theorem \ref{theo:caracdemaxasint}.
\begin{corollary}\label{cor:nomasquec*}
No estimator can guarantee a confidence level greater than $c^*$, given
by (\ref{eq:c*fud}), in an interval $(0,p_2)$.
\end{corollary}

According to Corollary \ref{cor:nomasquec*}, $c^*$ is an upper
bound on the confidence that can be guaranteed either asymptotically or
globally. It remains to be seen
% Vale esta frase
if there exists some estimator that can actually guarantee the
confidence level $c^*$. If it exists, that estimator will be optimum
from the point of view of guaranteed confidence. A related question is
if one such optimum estimator can be found within the restricted class
of non-randomized estimators. As will be shown, the answer to both
questions turns out to be affirmative for all values of $\mu_1$ and
$\mu_2$ in the case of asymptotic guarantee and for certain values of
$\mu_1$ and $\mu_2$ in the case of global guarantee.

Consider a non-randomized estimator $\hat p= g(N)$ of the form
%
%e3.1 ###
\begin{equation}
\label{eq:estimconsidgen}
g(n)= \frac{\Omega}{n+d},
\end{equation}
where $\Omega$ and $d$ are parameters (with $\Omega= \lim_{n
\rightarrow\infty} ng(n)$). This is a generalization
of the estimators (\ref{eq:hatpN-1n-1})--(\ref{eq:hatpNn}).
% *!*
Note that (\ref{eq:estimconsidgen}) has an asymptotic confidence
$\bar c$ given by (\ref{eq:barc}) and, for $\Omega= \Omega^*$, it
achieves the maximum $\liminf_{p\rightarrow0} c(p)$ that any
estimator can have, according to Theorem \ref{theo:caracdemaxasint}.

Under a mild condition on $d$, the estimator given by (\ref{eq:estimconsidgen}) with $\Omega= \Omega^*$ can be shown to asymptotically
guarantee the confidence $c^*$ for any $\mu_1$, $\mu_2>0$, as
established by Theorem \ref{theo:confgarasint} below.

\begin{theorem}
\label{theo:confgarasint}
The non-randomized estimator
%
%e3.2 ###
\begin{equation}
\label{eq:hatpOmegaoptn+d}
\hat p= \frac{\Omega^*}{N+d}
\end{equation}
asymptotically guarantees the optimum confidence $c^*$ given in (\ref
{eq:c*fud}) if
%
%e3.3 ###
\begin{equation}
\label{eq:dmayorque}
%d >  \{ -\nnum+ (\fud+1)/(\fud-1)  \} / 2,
d > \frac1 2  \biggl( -r+ \frac{M+1}{M-1}  \biggr),
\end{equation}
where $M$ is defined by (\ref{eq:fud}).
\end{theorem}

An estimator of this form can also \textit{globally} guarantee the
confidence level $c^*$, provided that $\mu_1$, $\mu_2$ are not too
small, as discussed in the following.

For $d \in\mathbb Z$, the estimator defined by (\ref{eq:estimconsidgen}) lends itself
% Vale
to an analysis similar to that carried out by \citeauthor{Mendo06} (\citeyear{Mendo06,Mendo08a}) for the particular case (\ref{eq:hatpN-1n}).
%(see the proof of lemma 1 of \citet{Mendo08a})
This allows the derivation of sufficient conditions on $\mu_1$, $\mu_2
$ which ensure that $c(p) \geq\bar c$ for all $p \in(0,1)$. The least
restrictive conditions are obtained for $d=1$ and are given in the
following proposition, which, particularized to $\Omega=\Omega^*$,
will yield the desired result on globally guaranteeing the optimum confidence.

\begin{proposition}
\label{prop:confgargen}
The confidence of the non-randomized estimator
%
%e3.4 ###
\begin{equation}
\label{eq:hatpOmegan+1}
\hat p= \frac{\Omega}{N+1}
\end{equation}
exceeds its asymptotic value $\bar c$, given by (\ref{eq:barc}), for
all $p \in(0,1)$ if
%
%e3.5 ###
\begin{equation}
\label{eq:condfunofdosgen}
\mu_1\geq\frac{\Omega}{r- \sqrt{r}}, \qquad\mu_2\geq
\frac{r+ \sqrt r+ 1}{\Omega}.
\end{equation}
\end{proposition}

Particularizing to $\Omega=\Omega^*$, Proposition \ref{prop:confgargen} establishes that there exists a non-randomized estimator that
can globally guarantee the optimum confidence $c^*$ for certain values
of $\mu_1$, $\mu_2$. It turns out that for $\Omega=\Omega^*$ and
$r$ given, one of the two inequalities in (\ref{eq:condfunofdosgen}) implies the other. Thus, for each $r$, only one of the
inequalities needs to be considered in order to determine the allowed
range for $\mu_1$, $\mu_2$. The result is stated in the following theorem.

\begin{theorem}
\label{theo:confgaropt}
The non-randomized estimator
%
%e3.6 ###
\begin{equation}
\label{eq:hatpOmegaoptn+1}
\hat p= \frac{\Omega^*}{N+1},
\end{equation}
with $\Omega^*$ as in (\ref{eq:Omega*}), globally guarantees the
optimum confidence $c^*$ given by (\ref{eq:c*fud}) if either of the
following conditions is satisfied:
%
%e3.8 ###
%e3.7 ###
\begin{eqnarray}
\label{eq:condoptfud34}
%(\fud-1)/\log\fud& \geq(\nnum+ \sqrt\nnum)/(\nnum-1), \quad\nnum
\frac{M-1}{\log M} &\geq&\frac{r+ \sqrt r}{r
-1} \qquad\mbox{for } r\in\{3,4\}; \\
\label{eq:condoptfud5omas}
\frac{M\log M}{M-1} &\geq&\frac{r+ \sqrt r+
1}{r} \qquad\mbox{for } r\geq5,
\end{eqnarray}
where $M$ is defined by (\ref{eq:fud}). These conditions can be
jointly expressed as $M\geq h(r)$, where $h$ is an increasing function.
\end{theorem}

Given $r$, consider the region of $(\mu_1,\mu_2)$ values that
satisfy the appropriate condition (\ref{eq:condoptfud34}) or
(\ref{eq:condoptfud5omas}), with $\mu_1,\mu_2>1$. From
Theorem \ref{theo:confgaropt}, the boundary of this region is a
continuous, decreasing, concave curve in $(1,\infty) \times(1,\infty
)$, namely, the curve determined by the equation $\mu_1\mu_2=h(r
)$. The region in question is the union of this curve and the portion
of the plane lying above and to the right. Furthermore, the region for
$r'>r$ contains that for $r$.

According to the preceding results, the problem of optimum estimation
of $p$, in the sense of globally guaranteeing the maximum possible
confidence, is solved by the non-randomized estimator
(\ref{eq:hatpOmegaoptn+1}) for $\mu_1$, $\mu_2$ satisfying
(\ref{eq:condoptfud34}) or (\ref{eq:condoptfud5omas}). Equivalently,
this estimator is minimax with respect to the risk defined by the loss
function given in (\ref{eq:lossfunction}) (maximin with respect to
confidence).

The fact that $c^*$ depends on $\mu_1$, $\mu_2$ only through $M$
gives rise to another interpretation of the result in Theorem \ref
{theo:confgaropt}. For $r$ and $M$ given, consider the
problem of finding, among all interval estimators of $p$ with a ratio
$M$ between their end-points, that which maximizes the globally
guaranteed confidence. The solution, if $M$ satisfies (\ref{eq:condoptfud34}) or (\ref{eq:condoptfud5omas}), is
%
%e3.9 ###
\begin{equation}
\label{eq:intest}
 \biggl[
\frac{r\log M}{(M-1)(N+1)},
\frac{rM\log M}{(M-1)(N+1)}
 \biggr]
\end{equation}
and the resulting maximum is $c^*$, as expressed by (\ref{eq:c*fud}). Equivalently, given $r$ and a prescribed confidence $c_0$,
if a value for $M$ is computed such that (\ref{eq:c*fud}) holds
with $c^*=c_0$ and if it satisfies (\ref{eq:condoptfud34}) or
(\ref{eq:condoptfud5omas}), then the interval estimator (\ref
{eq:intest}) minimizes the ratio between interval end-points subject
to a globally guaranteed confidence level $c_0$. Observe that it is
meaningful to prescribe, or minimize, the ratio of the interval
end-points, rather than their difference, since a given value for the
latter might be either unacceptably high or unnecessarily small,
depending on the unknown~$p$, whereas the ratio has a definite meaning,
regardless of $p$.

According to Proposition \ref{prop:confgargen}, conditions (\ref
{eq:condfunofdosgen}) are sufficient; however, they may not be
necessary. The same applies to (\ref{eq:condoptfud34}) and (\ref
{eq:condoptfud5omas}). Determining the most general conditions
which assure optimality of (\ref{eq:hatpOmegaoptn+1}) is a
difficult problem.\footnote{Although $c(p)$ (or a lower bound thereof)
% vale "thereof"
can be expressed in terms of the Gauss hypergeometric function
$_2F_1(a,b;c;t)$ (see the proof of Theorem \ref{theo:confgarasint}), standard algorithms for evaluation of hypergeometric sums
\citep{Petkovsek96} are not directly applicable because of the
existing dependence between $b$ and $t$.} Nevertheless, as will be
illustrated, the sufficient conditions (\ref{eq:condoptfud34}) or
(\ref{eq:condoptfud5omas}) cover most cases of interest.

%s4 ###
\section{Discussion}
\label{parte:disc}

\begin{figure}
\begin{tabular}{c}

\includegraphics{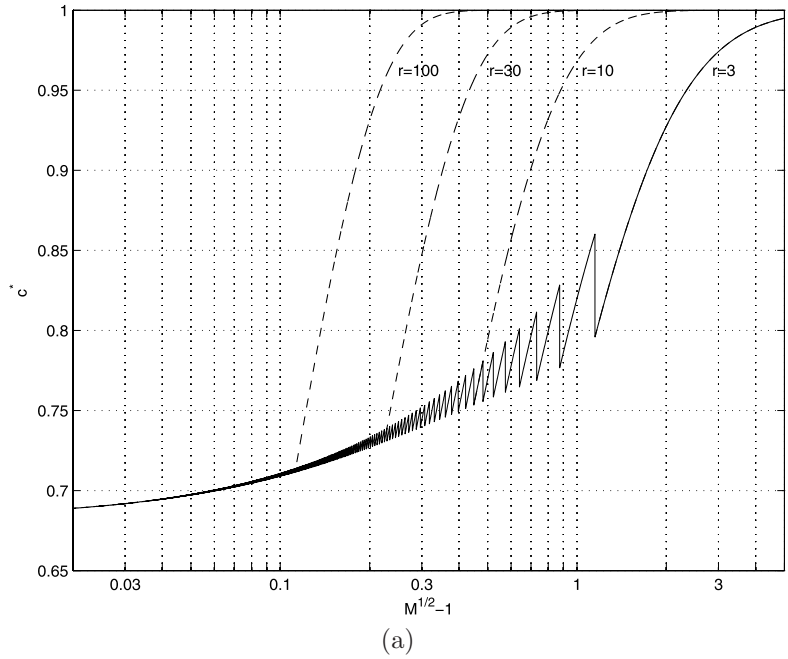}
\\
(a)\\

\includegraphics{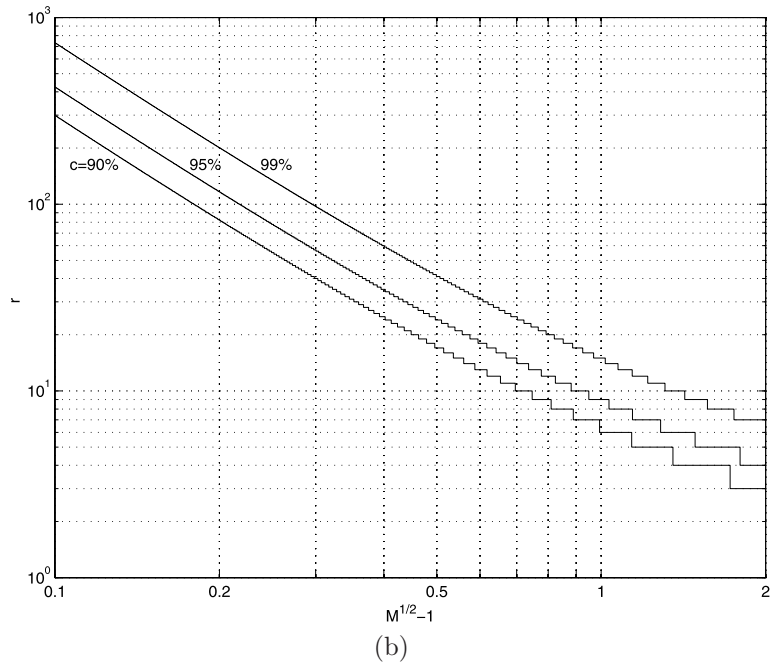}
\\
(b)
\end{tabular}
\caption{Behavior of optimum guaranteed confidence as a function of
$\sqrt M-1$: (a) $c^*$ (dashed) and minimum $c^*$ (solid); (b)
minimum $r$ that guarantees a confidence level.}\label{fig1}
\end{figure}

Figure~\ref{fig1}(a) depicts the relationship
% vale
between $c^*$, $M$ and $r$, for $M$ satisfying (\ref{eq:condoptfud34}) or (\ref{eq:condoptfud5omas}). The guaranteed
confidence $c^*$ is represented, for convenience, as a function of
$\sqrt M-1$, each dashed curve corresponding to a different $r
$. The figure also represents, with solid line, the minimum $c^*$ that
fulfills inequalities (\ref{eq:condoptfud34}) or (\ref{eq:condoptfud5omas}); this corresponds to the lowest $r$ for which the
applicable inequality holds. Figure~\ref{fig1}(b) shows this
minimum $r$ as a function of $\sqrt M-1$, with $c^*$ as a
parameter. From Corollary \ref{cor:nomasquec*}, this figure also
has a more general interpretation as the minimum~$r$ that is
required in order to guarantee (either globally or asymptotically) a
desired confidence level using any estimator based on inverse binomial sampling.

%f1 ###

Given $\mu_1$, $\mu_2$ and $c$, if the point $(\sqrt M-1,c)$ with
$M=\mu_1\mu_2$ lies in the region above the solid curve in
Figure~\ref{fig1}(a), then there exist values of $r$ for
which the estimator (\ref{eq:hatpOmegaoptn+1}) globally
guarantees the confidence level $c$ for the relative interval defined
by $\mu_1$ and $\mu_2$; the minimum such $r$ is displayed in
Figure~\ref{fig1}(b). As mentioned earlier, the region
referred to in Figure~\ref{fig1}(a) (or, equivalently,
conditions (\ref{eq:condoptfud34}) and (\ref{eq:condoptfud5omas})) covers most cases of interest. For example, any confidence
greater than $85\%$ can be globally guaranteed for any relative
interval with $\sqrt{M}-1 \leq1.108$, that is, such that $\mu_1
\mu_2\leq4.443$.

An important subclass of relative intervals is that for which $\mu_1
=\mu_2=1+m$,  $m>0$. An interval of this form corresponds to the
requirement that $\hat p$ and $p$ do not deviate from each other by a
factor greater than $1+m$; the parameter $m$ is thus interpreted as a
relative error margin. The guaranteed confidence and required $r$
in this case can be read directly from Figures \ref{fig1}(a)
and~\ref{fig1}(b), as $\sqrt{M}-1=m$. This particular
case is analyzed by \citeauthor{Mendo06} (\citeyear{Mendo06,Mendo08a}) for the estimator (\ref
{eq:hatpN-1n}). Comparing Figure~1 of \citet{Mendo08a} with
Figure~\ref{fig1}(a), the individual (dashed) curves in the
latter are seen to be above those in the former (this is most
noticeable for small $r$), in accordance with the fact that the
estimator (\ref{eq:hatpOmegaoptn+1}) is optimum. This yields a
reduction in the error margin $m$ for a given $r$ and a desired
guaranteed confidence. For example, taking $r=10$, the estimator
(\ref{eq:hatpN-1n}) guarantees a confidence level of $95\%$ for
$m=0.9074$, whereas (\ref{eq:hatpOmegaoptn+1}) guarantees the
same confidence for $m=0.8808$. Furthermore, the latter value is the
smallest $m$ for which a confidence level of $95\%$ can be guaranteed
by any estimator with the~$r$ in question.

As another manifestation of the optimum character of the estimator
(\ref{eq:hatpOmegaoptn+1}), the curves in Figure~\ref{fig1}(b) for $\sqrt{M}-1=m$ are farther to the left than those in
Figure~5(a) of \citet{Mendo06}. Since, from (\ref{eq:f}),
$\mathrm{E}
[N]=r/p$, for certain combinations of $m$ and $c$, this
provides a reduction in average observation time to achieve a globally
guaranteed confidence $c$ for an error margin $m$. (The fact that the
reduction is obtained only for certain combinations of $m$ and $c$ is a
consequence of the discrete character of $r$.) Thus, for $m=50\%$,
the estimator (\ref{eq:hatpN-1n}) requires $r=18$ in order to
globally guarantee a $90\%$ confidence level, whereas $r=17$
suffices for the estimator (\ref{eq:hatpOmegaoptn+1});
furthermore, this is the lowest required $r$ that can be achieved
by any estimator.

Comparing the attainable region shown in Figure~\ref{fig1}(a),
that is, the region above the solid curve, with the corresponding
region in Figure~1 of \citet{Mendo08a}, they are seen to have similar
shape and size, except that for small $m$, the boundary curve is
slightly higher in Figure~\ref{fig1}(a). Thus, the
applicability of the estimator (\ref{eq:hatpOmegaoptn+1}) is
similar to that of (\ref{eq:hatpN-1n}) (while achieving better
performance).

Another interesting particularization is $\mu_1=1+m$, $\mu_2
=1/(1-m)$, $0<m<1$, which corresponds to requiring that the absolute
error $|\hat p-p|$ does not exceed a fraction $m$ of the true value~$p$. The results for this case can be compared with those of \citet
{Mendo08a} and \citet{Chen07env}. For $m=40\%$ with a globally
guaranteed confidence level of $90\%$, the estimator~(\ref{eq:hatpN-1n}) requires $r=17$ (\citet{Mendo08a}, Proposition 1). The
results of \citet{Chen07env}, Theorem 2, give a sufficient value of
$r=46$ for the estimator (\ref{eq:hatpNn}). On the other
hand, from Theorem \ref{theo:confgaropt}, the estimator (\ref{eq:hatpOmegaoptn+1}) only requires $r=16$; moreover, it is assured
that any other estimator requires at least this $r$ to guarantee
(either globally or asymptotically) the same confidence for the $m$
under consideration.

%s5 ###
\section{Proofs}
\label{parte:proofs}

The following notation is introduced for convenience:
%
%e5.2 ###
%e5.1 ###
\begin{eqnarray}
\label{eq:phi}
\phi(t) &=& \frac{t^{r-1} \exp(-t)}{(r-1)!}, \\
\label{eq:ndunonddos}
n_1&=&  \lceil\Omega/(p\mu_1) - d  \rceil, \qquad
n_2=  \lfloor\Omega\mu_2/p - d  \rfloor.
\end{eqnarray}

\begin{pf*}{Proof of Proposition \ref{prop:asintnoasint}}
Let $c_\mathrm i$ denote $\liminf_{p \rightarrow p_0} c(p)$ and assume
that $c_\mathrm i < c_0$. For $\epsilon= (c_0 - c_\mathrm i)/2$, the
definition of limit inferior implies that there exists $p_\epsilon\in
(p_1,p_2)$ such that $c(p_\epsilon) < c_\mathrm i + \epsilon< c_0$
and thus the estimator does not guarantee the confidence level $c_0$ in
$(p_1,p_2)$. This establishes the result.
\end{pf*}

\begin{lemma}\label{lemma:dergammax}
For all $t>0$, the function $\phi$ defined in (\ref{eq:phi}) satisfies
%
%e5.3 ###
\begin{equation}
\label{eq:cotaendergamma}
0 < \phi(t) \leq\frac{(r-1)^{r-1} \exp(-r+1)}{(r
-1)!} < 1.
\end{equation}
\end{lemma}

\begin{pf}
The first inequality in (\ref{eq:cotaendergamma}) is obvious. Let
$
Q_r= (r-1)^{r-1} \exp(-r+1) /(r-1)!.
$
Maximizing $\phi(t)$ with respect to $t$, it is seen that $\phi(t)
\leq Q_r$. Since $Q_{r+1}/Q_r= (1 + 1/(r-1))^{r
-1} \exp(-1) < 1$ and $Q_3 = 2\exp(-2) < 1$, it follows that $Q_r
< 1$ for all $r\geq3$.
\end{pf}

\begin{pf*}{Proof of Proposition \ref{prop:limdeconf}}
Consider $g\in\mathcal{F}_{\mathrm{n}}$ and let $\Omega= \lim
_{n\rightarrow
\infty} ng(n)$. Given $\epsilon>0$, there exists
$n_\epsilon$ such that $|ng(n)-\Omega|<\epsilon$
for all $n\geq n_\epsilon$, that is,
$
(\Omega-\epsilon)/n< g(n) < (\Omega+\epsilon)/n.
$
Using this, the confidence $c(p) = P[p/\mu_2\leq g(N)
\leq p \mu_1]$ can be bounded for $p \leq(\Omega-\epsilon)/(\mu_1
n_\epsilon)$ as
%
%e5.5 ###
%e5.4 ###
\begin{eqnarray}
\label{eq:cotainfcpFeln}
c(p) & \geq & P [  \lceil(\Omega+\epsilon)/(p\mu_1)
 \rceil\leq N\leq \lfloor(\Omega-\epsilon)\mu_2
/p  \rfloor ], \\
\label{eq:cotasupcpFeln}
c(p) & \leq & P [  \lceil(\Omega-\epsilon)/(p\mu_1)
 \rceil\leq N\leq \lfloor(\Omega+\epsilon)\mu_2
/p  \rfloor ].
\end{eqnarray}

Let $b(n,p;i)$ denote the binomial probability function with
parameters $n$ and $p$ evaluated at~$i$. From the relationship
between binomial and negative binomial distributions, (\ref{eq:cotainfcpFeln}) is written as
%
%e5.6 ###
\begin{equation}
\label{eq:cotainfcpFeln2}
c(p) \geq
% \rfloor ] -
% \rceil- 1  ] \\ =
\sum_{i=0}^{r-1} b \bigl(  \lceil(\Omega+\epsilon
)/(p\mu_1)  \rceil-1,p;i  \bigr) -
\sum_{i=0}^{r-1} b \bigl(  \lfloor(\Omega-\epsilon
)\mu_2/p  \rfloor,p;i  \bigr).
\end{equation}
According to the Poisson theorem (\citet{Papoulis02}, page 113), the
right-hand side of~(\ref{eq:cotainfcpFeln2}) converges to
$\gamma(r,(\Omega-\epsilon)\mu_2) - \gamma(r,(\Omega
+\epsilon)/\mu_1)$ as $p \rightarrow0$. This implies that
%
%e5.7 ###
\begin{equation}
\label{eq:cotainfcpFeln3}
\liminf_{p \rightarrow0} c(p) \geq\gamma\bigl(r,(\Omega-\epsilon
)\mu_2\bigr) - \gamma\bigl(r,(\Omega+\epsilon)/\mu_1\bigr).
\end{equation}
Since $\partial\gamma(r,t) / \partial t = \phi(t)$, Lemma \ref
{lemma:dergammax} establishes that $0 < \partial\gamma(r,t) /
\partial t < 1$. Using this, (\ref{eq:cotainfcpFeln3}) yields
$
\liminf_{p \rightarrow0} c(p) \geq\bar c - (\mu_2+1/\mu_1)\epsilon.
$
Taking into account that this holds for all $\epsilon>0$, it follows that
$
\liminf_{p \rightarrow0} c(p) \geq\bar c.
$
Applying similar arguments to (\ref{eq:cotasupcpFeln}) gives
$
\limsup_{p \rightarrow0} c(p) \leq\bar c.
$
Thus, $\lim_{p \rightarrow0} c(p)$ exists and equals $\bar c$.
\end{pf*}

\begin{lemma}[(\citet{Abramowitz70}, inequality (4.1.33))] %
%parece que hacen falta dos pares de llaves extra! / It appears two
%extra pairs of braces are required!
\label{lemma:ln}
$
t/(1+t) \leq\log(1+t) \leq t
$
for $t>-1$, with equality if and only if $t=0$.
\end{lemma}

\begin{lemma}
\label{lemma:convunif}
Let $(p_k)$ be a positive sequence which converges to $0$. For any $\nu
_1$, $\nu_2$ with $\nu_2 > \nu_1 > 0$, the sequence of functions
$(\phi_k)$ defined as
%
%e5.8 ###
\begin{equation}
\label{eq:phik}
\phi_k(\nu) = \frac{(1-p_k)^{\nu/p_k-r}}{(r-1)!} \prod
_{i=1}^{r-1} (\nu-ip_k)
\end{equation}
converges uniformly to $\phi(\nu)$, given by (\ref{eq:phi}), in the
interval $[\nu_1, \nu_2]$.
\end{lemma}

\begin{pf}
Since the sequence $(p_k)$ is positive and converges to $0$, there
exists $k_1$ such that $p_k<\min\{\nu_1/(r-1),1\}$ for $k \geq
k_1$. Thus, $\phi_k(\nu)>0$ for $\nu\in[\nu_1,\nu_2]$, $k \geq
k_1$ and, therefore,
%
%e5.9 ###
\begin{equation}
\label{eq:phiphikminmax1}
|\phi_k(\nu)-\phi(\nu)| = \min\{\phi_k(\nu),\phi(\nu)\}
\biggl(\max \biggl\{ \frac{\phi_k(\nu)}{\phi(\nu)}, \frac{\phi(\nu
)}{\phi_k(\nu)}  \biggr\}-1 \biggr).
\end{equation}
Lemma \ref{lemma:dergammax} implies that
$
\min\{\phi_k(\nu),\phi(\nu)\} \leq1
$.
On the other hand,
\[
\max \biggl\{ \frac{\phi_k(\nu)}{\phi(\nu)}, \frac{\phi(\nu
)}{\phi_k(\nu)}  \biggr\} = \exp \biggl|\log\frac{\phi_k(\nu
)}{\phi(\nu)}  \biggr|.
\]
Substituting into (\ref{eq:phiphikminmax1}),
we have
%
%e5.10 ###
\begin{equation}
\label{eq:phiphikminmax2}
|\phi_k(\nu)-\phi(\nu)| \leq\exp \biggl|\log\frac{\phi_k(\nu
)}{\phi(\nu)}  \biggr| - 1.
\end{equation}

From (\ref{eq:phi}) and (\ref{eq:phik}), it follows that, assuming
$\nu\in[\nu_1,\nu_2]$ and $k \geq k_1$,
%
%e5.11 ###
\begin{equation}
\label{eq:lnrelphi}
 \biggl| \log\frac{\phi_k(\nu)}{\phi(\nu)}  \biggr| \leq \Biggl|
\sum_{i=1}^{r-1} \log \biggl(1-\frac{ip_k}{\nu}  \biggr)
\Biggr| +  \biggl|  \biggl( \frac{\nu}{p_k} - r \biggr) \log(1-p_k) +
\nu \biggr|.
\end{equation}
Using Lemma \ref{lemma:ln}, the first term in the right-hand side of
(\ref{eq:lnrelphi}) is bounded as
follows:% pongo "<", no "<=", porque p_k es > 0
%
%e5.12 ###
\begin{equation}
\label{eq:cotaconvunif1}
 \Biggl| \sum_{i=1}^{r-1} \log \biggl(1- \frac{ip_k}{\nu}
\biggr)  \Biggr| = -\sum_{i=1}^{r-1} \log \biggl(1-\frac{ip_k}{\nu}
 \biggr) < \sum_{i=1}^{r-1} \frac{ip_k}{\nu-ip_k}.
\end{equation}
As for the second term, since $\nu\in[\nu_1,\nu_2]$, there exists
$k_2 \geq k_1$ such that $p_k< \min\{\nu_1/r,1\}$ for $k \geq
k_2$. Thus, $(\nu/p_k-r)\log(1-p_k)<0$, whereas $\nu> 0$. From
Lemma \ref{lemma:ln}, $-p_k/(1-p_k)<\log(1-p_k)<-p_k$. It follows that
%
%e5.13 ###
\begin{eqnarray}\label{eq:cotaconvunif3}
\hspace*{-5pt}\biggl|  \biggl( \frac{\nu}{p_k} - r \biggr) \log(1-p_k) + \nu\biggr|
&<& \max \biggl\{  \biggl| - \biggl( \frac{\nu}{p_k} - r
 \biggr) p_k + \nu \biggr|,  \biggl| - \biggl( \frac{\nu}{p_k} - r
 \biggr) \frac{p_k}{1-p_k} + \nu \biggr|  \biggr\}\qquad\hspace*{8pt}
 \nonumber\\[-8pt]\\[-8pt]
&=& \max \biggl\{ r, \frac{|r-\nu|}{1-p_k}  \biggr\}
p_k < \max\{r,\nu_2\} \frac{p_k}{1-p_k}.\nonumber
\end{eqnarray}
Combining (\ref{eq:phiphikminmax2})--(\ref{eq:cotaconvunif3}) yields the following bound, valid for all $\nu\in[\nu_1,\nu_2]$,
$k \geq k_2$:
%
%e5.14 ###
\begin{equation}
\label{eq:cotaconvunif5}
|\phi_k(\nu) - \phi(\nu)| < \exp \Biggl[  \Biggl( \sum_{i=1}^{r
-1} \frac{i}{\nu-ip_k} + \frac{\max\{r,\nu_2\}}{1-p_k}
\Biggr) p_k  \Biggr] - 1.
\end{equation}
The right-hand side of (\ref{eq:cotaconvunif5}) tends to $0$ as $k
\rightarrow\infty$. Therefore, $\sup_{\nu\in[\nu_1,\nu_2]} |\phi
_k(\nu) - \phi(\nu)|$ also tends to $0$, which establishes the result.
\end{pf}

\begin{pf*}{Proof of Theorem \ref{theo:caracdemaxasint}}
For $\alpha\in(0, 1]$, let the sequence $(p_k)$, $k \in\mathbb N$ be
defined as $p_k = \alpha M^{-k}$, with $M$ given by (\ref{eq:fud}),
and let the sequence of intervals $(I_k)$ be defined as $I_k =
(p_k/\mu_2, p_k\mu_1]$. Consider also the sequence $(f_k)$, where
$f_k$ is the probability function of $N$ with parameters $r
$ and $p_k$. From (\ref{eq:f}) and (\ref{eq:phik}),
%
%e5.15 ###
\begin{equation}
\label{eq:fkphik}
f_k(n) = p_k \phi_k(np_k).
\end{equation}

To facilitate the development, it is convenient to first analyze
non-randomized estimators, and then generalize to randomized estimators.
Given a non-randomized estimator specified by $g\in\mathcal{F}$, let the
sequence of sets $(S_k)$ be defined such that $n\in S_k$ if and
only if $g(n) \in I_k$. Since the intervals $I_k$ are
disjoint for different $k$, the sets $S_k$ are also disjoint. Let the
function $\sigma$ be defined such that $\sigma(n)=k$ if and only
if $n\in S_k$, with $\sigma(n)=0$ (or an arbitrary negative
value) if $n\notin S_k$ for all $k \in\mathbb N$. Thus, $\sigma$
gives the index $k$ of the interval $I_k$ that $g$ associates to
each $n$, if any. The function $\sigma$ is determined by $g
$; and, for a given $n_0$, $\sigma(n_0)$ can be modified
without affecting the rest of values $\sigma(n)$, $n\neq
n_0$ by adequately choosing $g(n_0)$.

For the considered non-randomized estimator, let $c_k$ denote the
probability that $\hat p$ lies in $I_k$ when $p = p_k$, i.e. $c_k =
\sum_{n\in S_k} f_k(n)$. Further, let
$
C(h,H) = \sum_{k=h}^{h+H-1} c_k
$
for $h$, $H \in\mathbb N $ arbitrary. Defining
\[
s(n) =
\cases{
f_{\sigma(n)}(n), &\quad if $ n\in S_h \cup\cdots
\cup S_{h+H-1}$, \cr
0, &\quad otherwise,
}
\]
the sum $C(h,H)$ can be expressed as $\sum_{n= r}^\infty
s(n)$. It follows that
%
%e5.16 ###
\begin{equation}
\label{eq:CHcota}
C(h,H) \leq\sum_{n=r}^\infty  \max_{k \in\{h,\ldots
,h+H-1\}} f_k(n),
\end{equation}
which holds with equality if the estimator satisfies
%
%e5.17 ###
\begin{equation}
\label{eq:reglamaxCH}
\sigma(n) = \arg\max_{k \in\{h,\ldots,h+H-1\}} f_k(n)
\qquad\mbox{for all } n\in\mathbb N, n\geq r,
\end{equation}
where, if the maximum is reached at more than one index $k$, the $\arg
\max$ function is arbitrarily defined to give the lowest such index.

As for randomized estimators, consider an arbitrary function $G
\in\mathcal{F}_{\mathcal{R}}$, that for each $n$ specifies $\Pi_n$, the
distribution function of $\hat p$ conditioned on $N= n$.
Let $I^\mathrm c = (-\infty, p_{h+H-1}/\mu_2] \cup(p_h\mu_1, \infty
)$. Conditioned on $N= n$, let $\pi_{n,k}$ and $\pi
_n^\mathrm c$ respectively denote the probabilities that $\hat p
$ is in $I_k$ and in $I^\mathrm c$. Obviously,
%
%e5.18 ###
\begin{equation}
\label{eq:pisuman1}
\pi_{n,h}+\cdots+\pi_{n,h+H-1} + \pi_n^\mathrm c = 1.
\end{equation}
For $N=n$, the numbers $\pi_{n,h}, \ldots, \pi
_{n,h+H-1}, \pi_n^\mathrm c$ indicate how the conditional
probability associated to all possible values of $\hat p$ (that is,
$1$ in total) is divided among $I_h$, \ldots, $I_{h+H-1}$, $I^\mathrm
c$. For a given $n_0$, any combination of values $\pi_{n
_0,h}, \ldots, \pi_{n_0,h+H-1}, \pi_{n_0}^\mathrm c$
allowed by (\ref{eq:pisuman1}) can be realized, without affecting
other values $\pi_{n,k}$, $\pi_{n}^\mathrm c$ for $n
\neq n_0$, by adequately choosing the distribution function $\Pi
_{n_0}$. Defining $c_k$ and $C(h,H)$ as in the non-randomized case,
the former is expressed as
%
%e5.19 ###
\begin{equation}
\label{eq:ckrand}
c_k = \sum_{n=r}^\infty\pi_{n,k} f_k(n).
\end{equation}
Since all terms in (\ref{eq:ckrand}) are positive, the series
converges absolutely, and thus $C(h,H)$ can be written as
\[
C(h,H) = \sum_{k=h}^{h+H-1} c_k = \sum_{n=r}^\infty\sum
_{k=h}^{h+H-1} \pi_{n,k} f_k(n).
\]
It is evident from this expression that $C(h,H)$ is maximized if, for
each $n$, the values $\pi_{n,k}$, $\pi_n^\mathrm c$
are chosen as $\pi_{n,l}=1$ for $l = \arg\max_k f_k(n)$,
$\pi_{n,k}=0$ for $k \neq l$, $\pi_n^\mathrm c = 0$; and
the resulting maximum coincides with the right-hand side of (\ref{eq:CHcota}). Thus the inequality (\ref{eq:CHcota}), initially derived
for non-randomized estimators, also holds for the general class of
randomized estimators. This implies that for any randomized estimator
there exists a non-randomized estimator that attains the same or greater
$C(h,H)$.

Let the function $\bar g\in\mathcal{F}$ be defined as $\bar g(n)
= \Omega^*/n$, and consider the non-randomized estimator specified
by this function. The sets $S_k$ and function $\sigma$ associated to
this estimator will be denoted as $\bar S_k$ and $\bar\sigma$
respectively. It is seen that
%
%e5.20 ###
\begin{equation}
\label{eq:tildeSk}
\bar S_k = \bigl[\Omega^*/(\mu_1p_k), \Omega^*\mu_2/p_k\bigr) \cap\{r,
r+1, r+2,\ldots\}.
\end{equation}
Let the sequence of functions $(\bar f_k)$ be defined as
$
\bar f_k(n) = p_k \phi(np_k).
$
It is readily seen that the equation $\bar f_k(t) = \bar f_{k+1}(t)$ has
only one solution, given as $t=\Omega^*\mu_2/p_k$. Similarly, $\bar
f_{k-1}(t) = \bar f_{k}(t)$ has the solution $t=\Omega^*/(\mu_1p_k)$.
Taking into account that the functions $\bar f_k$ are unimodal, this
implies that $\bar f_k(n) \geq\bar f_{k-1}(n)$ and $\bar
f_k(n) \geq\bar f_{k+1}(n)$ for $n\in\bar S_k$. An
analogous argument shows that $\bar f_k(n) > \bar f_{k-i}(n)$
and $\bar f_k(n) > \bar f_{k+i}(n)$ for $n\in\bar S_k$ and
$i \geq2$. Therefore, the function $\bar\sigma$ associated to $\bar
g$ satisfies a modified version of (\ref{eq:reglamaxCH}) in
which $f_k$ is replaced by $\bar f_k$, that is, $\phi_k$ in (\ref{eq:fkphik}) is replaced by its limit $\phi$.
% (see upper and middle parts of Figure~\ref{fig:demo}) <Poner slo si
%pongo la figura>.
In the following, using the fact that the difference between $\phi_k$
and $\phi$ is small for large $k$, the estimator defined by $\bar g
$ will be used to derive from (\ref{eq:CHcota}) a more explicit
upper bound on $C(h,H)$. Since the sum $C(h,H)$ attained by any
estimator is equalled or exceeded by some non-randomized estimator, it
will be sufficient to restrict to the class of non-randomized estimators.

Consider an arbitrary non-randomized estimator defined by $g\in\mathcal{F}
$ with its corresponding function $\sigma$. Differentiating (\ref{eq:f}) with respect to $p$, it is seen that $\partial f(n) /
\partial p$ is positive for $p < r/n$ and negative for $p >
r/n$. According to (\ref{eq:tildeSk}),
%
%e5.21 ###
\begin{equation}
\label{eq:cotasndenpsigma0}
\Omega^* / (\mu_1n) \leq p_k < \Omega^*\mu_2/ n\qquad
\mbox{for }n\in\bar S_k.
\end{equation}
Using Lemma \ref{lemma:ln}, it stems from (\ref{eq:Omega*funofdos}) that
%
%e5.22 ###
\begin{equation}
\label{eq:funoOmega*fdos}
\Omega^*/\mu_1< r< \Omega^*\mu_2.
\end{equation}
The first inequality in (\ref{eq:cotasndenpsigma0}) and the
second in (\ref{eq:funoOmega*fdos}) imply that $p_{k-1} = Mp_k
\geq\Omega^*\mu_2/n> r/n$ for $n\in\bar S_k$.
Thus $\partial f(n) / \partial p < 0$ for $p \geq p_{k-1}$. This
implies that $f_{k+i}(n) < f_{k-1}(n)$ for $n\in\bar
S_k$ and $i \leq-2$. It follows that if $\sigma(n_0) < \bar
\sigma(n_0)-1$ for a given $n_0$, the sum $C(h,H)$ could be
made larger by modifying the value $g(n_0)$ so as to attain
$\sigma(n_0) = \bar\sigma(n_0)-1$; unless $\sigma(n
_0)=h+H-1$, in which case modifying $g(n_0)$ cannot make
$C(h,H)$ larger. Analogously, the second inequality in (\ref{eq:cotasndenpsigma0}) and the first in (\ref{eq:funoOmega*fdos}) yield
%Likewise, the second inequality in \eqref{eq: cotas nden p sigma 0},
%together with the first in \eqref{eq: funo Omega * fdos}, implies that
%$p_{k+1} = p_k/\fud< \Omega^*/(\funo\nden) < \nnum/\nden$ for $\nden
$f_{k+i}(n) < f_{k+1}(n)$ for $n\in\bar S_k$ and $i \geq
2$. Therefore, if $\sigma(n_0) > \bar\sigma(n_0)+1$ for
some $n_0$, $C(h,H)$ could be made larger by modifying $g
(n_0)$, unless $\sigma(n_0)=h$.

%Figure~\ref{fig:demo} illustrates the arguments in the preceding
%paragraph. The set $\ea S_k$ as well as the functions $f_{k+i}$ and $
%different horizontal axes for clarity. Observe that the curves
%$f_{k+i}$ are slightly modified versions of $\ea f_{k+i}$. For $\nden
%f_{k+i}(\nden)$ for all $i$, $f_k(\nden)$ does not necessarily exceed
%$f_{k-1}(\nden)$ (see shaded area) or $f_{k+1}(\nden)$; however, it
%holds that $f_{k-1}(\nden) > f_{k-2}(\nden)$ and $f_{k+1}(\nden) >
%f_{k+2}(\nden)$.

%Illustration of inequalities $f_{k-i}(\nden) < f_{k-i+1}(\nden)$ and
%$f_{k+i}(\nden) < f_{k+i-1}(\nden)$ for $\nden\in\ea S_k$, $i \geq
%2$}%

According to the above, in order to obtain an upper bound on $C(h,H)$,
it suffices to consider non-randomized estimators such that $\sigma
(n) \in\{\bar\sigma(n)-1, \bar\sigma(n), \bar\sigma
(n)+1\}$ for $n\in\bar S_h \cup\cdots\cup\bar S_{h+H-1}$,
$\sigma(n) = h$ for $n< \lceil\Omega^*/(\mu_1p_h) \rceil
$ and $\sigma(n) = h+H-1$ for $n\geq\lceil\Omega^*\mu_2
/p_{h+H-1} \rceil$. Thus, from (\ref{eq:CHcota}),
%
%e5.23 ###
\begin{eqnarray}
\label{eq:cotaChHcon2}
C(h,H) & \leq&\sum_{k=h}^{h+H-1} \sum_{n\in\bar S_k} \max_{i
\in\{-1,0,1\}} f_{k+i}(n)
+ \sum_{n= r}^{\lceil\Omega^*/(\mu_1p_h) \rceil-1}
f_h(n)\nonumber\\
&&{}+ \sum_{n= \lceil\Omega^*\mu_2/p_{h+H-1} \rceil
}^\infty f_{h+H-1}(n)\\
&\leq&\sum_{k=h}^{h+H-1} \sum_{n\in\bar S_k} \max_{i \in\{
-1,0,1\}} f_{k+i}(n) + 2.\nonumber
\end{eqnarray}

Let $\nu_1 = \Omega^*/\mu_1$, $\nu_2 = \Omega^*\mu_2$, and
consider $\epsilon>0$. From Lemma \ref{lemma:convunif}, there
exists $k_1$ such that $|\phi_k(\nu)-\phi(\nu)| < \epsilon$ for
$\nu\in[\nu_1,\nu_2]$, $k \geq k_1$. In addition, (\ref{eq:tildeSk}) implies that $np_k \in[\nu_1, \nu_2]$ for $n\in\bar
S_k$. From these facts and (\ref{eq:fkphik}) it stems that
$
 | f_k(n)/p_k - \phi(np_k)  | < \epsilon
$
for $n\in\bar S_k$, $k \geq k_1$. Thus, for $k \geq k_1+1$
\begin{eqnarray}
\label{eq:cotamax1}
\max_{i \in\{-1,0,1\}} f_{k+i}(n)
& < &\max_{i \in\{-1,0,1\}} \bigl[ p_{k+i} \bigl(\phi(np_{k+i}) + \epsilon
\bigr) \bigr] \nonumber\\[-8pt]\\[-8pt]
& \leq&\max_{i \in\{-1,0,1\}} ( p_{k+i} \phi(np_{k+i})) +
\epsilon p_{k-1}.\nonumber
\end{eqnarray}
As previously shown, $\bar f_k(n) = p_k \phi(np_k)$ equals or
exceeds $\bar f_{k-1}(n)$ and $\bar f_{k+1}(n)$ for $n\in
\bar S_k$. Therefore, from (\ref{eq:cotamax1})
%
%e5.24 ###
\begin{equation}
\label{eq:cotamax2}
\max_{i \in\{-1,0,1\}} f_{k+i}(n) < p_{k} \phi(np_{k}) +
\epsilon p_{k-1}.
\end{equation}
Substituting (\ref{eq:cotamax2}) into (\ref{eq:cotaChHcon2}),
%
%e5.25 ###
\begin{equation}
\label{eq:CHmenor1}
C(h,H) \leq\sum_{k=h}^{h+H-1} \sum_{n\in\bar S_k} \bigl(p_k \phi
(np_k) + \epsilon p_{k-1}\bigr) + 2.
\end{equation}
The number of elements in the set $\bar S_k$ is less than $\Omega
^*(\mu_2-1/\mu_1)/p_k+1$, according to (\ref{eq:tildeSk}). From
(\ref{eq:Omega*funofdos}), this upper bound equals $r\log
M/p_k+1$. Using this in (\ref{eq:CHmenor1}), and taking into
account that $p_{k-1}<1$,
%
%e5.26 ###
\begin{equation}
\label{eq:CHmenor1bis}
C(h,H) < \sum_{k=h}^{h+H-1} \sum_{n\in\bar S_k} p_k \phi(n
p_k) + \epsilon H(rM\log M+1) + 2.
\end{equation}
The term $\sum_{n\in\bar S_k} p_k \phi(np_k)$ in (\ref{eq:CHmenor1bis}) converges to $\int_{\Omega^*/\mu_1}^{\Omega
^*\mu_2} \phi(\nu) \,\mathrm{d}\nu= c^*$ as $k \rightarrow\infty$.
Thus, for the considered $\epsilon$, there exists $k_2$ such that for
all $k \geq k_2$
%
%e5.27 ###
\begin{equation}
\label{eq:sumaintegral}
 \biggl| \sum_{n\in\bar S_k} p_k \phi(np_k) - c^*  \biggr| <
\epsilon.
\end{equation}
Consequently, defining $k_0 = \max\{k_1+1,k_2\}$, inequalities (\ref
{eq:CHmenor1bis}) and (\ref{eq:sumaintegral}) imply that for all
$h \geq k_0$, for all $H$, and for $\alpha\in(0,1]$,
%
%e5.28 ###
\begin{eqnarray}
\label{eq:CHmenor2}
C(h,H) &<& H(c^*+\epsilon) + \epsilon H (rM\log M+1) + 2
\nonumber
\\[-8pt]
\\[-8pt]
\nonumber
&=&
H (c^* + \epsilon P) + 2,
\end{eqnarray}
where $P = rM\log M+2>0$.

Since the minimum of a set cannot be larger than the average of the
set, from (\ref{eq:CHmenor2}) it follows that there is some $k_3
\in\{h,\ldots,h+H-1\}$ such that $c_{k_3} \leq C(h,H)/H < c^* +
\epsilon P + 2/H$.

Using the foregoing results, the bound $\liminf_{p \rightarrow0} c(p)
\leq c^*$ for an arbitrary estimator can be established by
contradiction. Assume that there is some estimator, defined by $G
\in\mathcal{F}_{\mathcal{R}}$, such that $\liminf_{p \rightarrow0}
c(p) = c^* + d$ with
$d>0$. This means that for any $\epsilon'>0$ there exists $p_0$ such
that $c(p) \geq c^* + d -\epsilon'$ for all $p \leq p_0$. Thus taking
$\epsilon' = d/3$, there is $p_0$ such that
%
%e5.29 ###
\begin{equation}
\label{eq:cgeqparacontrad}
c(p) \geq c^* + 2d/3 \qquad\mbox{for all } p \leq p_0.
\end{equation}
On the other hand, taking $\epsilon=d/(6P)$, $H = \lceil12/d \rceil
$, and with $\alpha\in(0,1]$ arbitrary, the result in the preceding
paragraph assures that for any $h$ not smaller than a certain $k_0$
(which depends on the considered $\epsilon$) there exists $k_3 \in\{
h,\ldots,h+H-1\}$ with
%
%e5.30 ###
\begin{equation}
\label{eq:ck1}
c_{k_3} < c^* + \epsilon P + 2/H \leq c^* + d/3.
\end{equation}

Let $h$ be selected such that
$
\label{eq:hgeqmin}
h \geq\max \{
k_0,
-\log p_0 / \log M
 \}.
$
For each $k = h,\ldots,h+H-1$, let $X_k$ denote the (possibly empty)
set of all points $x \in\mathbb R$ such that, for the considered
estimator and for $p=p_k$, the probability that $\hat p$ equals $x$
is at least $d/3$. The number of points in $X_k$ cannot exceed $\lfloor
3/d \rfloor$, for otherwise the sum of their probabilities would be
greater than $1$. The set $X_k$ is determined by $f_k$ and $G$
(or by $f_k$ and $g$, in the case of a non-randomized estimator
defined as $\hat p=g(N)$).

Let $\alpha\in(0,1]$ be chosen such that
%
%e5.31 ###
\begin{equation}
\label{eq:condalpha}
p_{k}/\mu_2\notin X_k \qquad\mbox{for all } k=h,\ldots,h+H-1.
\end{equation}
Such $\alpha$ necessarily exists because (\ref{eq:condalpha})
excludes only a finite number of possible values from $(0,1]$. This
choice of $\alpha$ assures that for $p=p_k$, $k =h,\ldots,h+H-1$, the
probability that $\hat p$ equals $p_k/\mu_2$ is smaller than $d/3$. Thus
$
c(p_k) < c_k + d/3
$
for $k = h,\ldots,h+H-1$, which, together with (\ref{eq:ck1}), gives
%
%e5.32 ###
\begin{equation}
\label{eq:cpk1}
c(p_{k_3}) < c^* + 2d/3.
\end{equation}
On the other hand, since $k_3 \geq h$, from the choice of $h$ it
follows that
\[
p_{k_3} = \alpha M^{-k_3} \leq\alpha M^{-h} \leq\alpha M
^{\log p_0 /\log M} = \alpha p_0 \leq p_0.
\]
This implies, according to (\ref{eq:cgeqparacontrad}), that
$c(p_{k_3}) \geq c^* + 2d/3$, in contradiction with (\ref{eq:cpk1}). Therefore $\liminf_{p \rightarrow0} c(p) \leq c^*$.

It has been shown that, for any estimator, $\liminf_{p \rightarrow0}
c(p)$ cannot exceed $c^*$. In addition, any non-randomized estimator
defined by a function $g\in\mathcal{F}_{\mathrm{n}}$ with $\lim
_{n\rightarrow
\infty} ng(n) = \Omega^*$ achieves $\liminf_{p
\rightarrow0} c(p) = \lim_{p \rightarrow0} c(p) = c^*$. Therefore,
$c^*$ is the maximum of $\liminf_{p \rightarrow0} c(p)$ over all
randomized or non-randomized estimators.
\end{pf*}

\begin{pf*}{Proof of Theorem \ref{theo:confgarasint}}
Let $I_p(z,w)$ denote the regularized incomplete beta function:
\begin{eqnarray*}
I_p(z,w) & = &\frac{1}{B(z,w)} \int_0^p t^{z-1}(1-t)^{w-1} \,\mathrm
{d}t,\\
B(z,w) & = &\int_0^1 t^{z-1}(1-t)^{w-1} \,\mathrm{d}t\\
& =& \frac{\Gamma
(z)\Gamma(w)}{\Gamma(z+w)}.
\end{eqnarray*}
From \citet{Abramowitz70}, equation (6.6.4),
$
P[N\leq n] = I_p(r,n-r+1)
$
for $n\in\mathbb N$,  $n\geq r$. The confidence for
the estimator (\ref{eq:hatpOmegaoptn+d}) can thus be written as
%
%e5.33 ###
\begin{equation}
\label{eq:confbetainc}
c(p) = P[n_1\leq N\leq n_2] = I_p(r,n_2
-r+1) - I_p(r,n_1-r),
\end{equation}
where $n_1$ and $n_2$ are given by (\ref{eq:ndunonddos})
with $\Omega=\Omega^*$. It is easy to show that $I_p(z,w)$ is an
increasing function of $w$ for $w \in\mathbb R$. Since $n_1<
\Omega^*/(p\mu_1) - d + 1$ and $n_2> \Omega^*\mu_2/p - d - 1$,
it follows from (\ref{eq:confbetainc}) that
%
%e5.34 ###
\begin{equation}
\label{eq:cc1c2opt}
c(p) > \tilde c_2(p) - \tilde c_1(p)
\end{equation}
with
$
\tilde c_1(p) = I_p  ( r, \Omega^*/(p\mu_1) - r- d + 1
 )
$ and
$
\tilde c_2(p) = I_p  ( r, \Omega^*\mu_2/p - r- d
 ).
$

Expressing $I_p(z,w)$ as (\citet{Abramowitz70}, equation (26.5.23))
\begin{eqnarray*}
I_p(z,w) & = & \frac{{}_2F_1(z,1-w;z+1;p) p^z}{B(z,w)z},\\
_2 F_1(a,b;c;t) & = & \sum_{j=0}^\infty\frac{ (a+j-1)^{(j)}
(b+j-1)^{(j)} t^j}{ (c+j-1)^{(j)} j!},
\end{eqnarray*}
the term $\tilde c_2(p)$ can be written as
%
%e5.36 ###
%e5.35 ###
\begin{eqnarray}
\label{eq:c2cota2antes}
\tilde c_2(p) & = & \frac{1}{B(r, \Omega^*\mu_2/p-r-d)}
\frac{p^r}{r}     {}_2F_1 (r, r+d+1-\Omega
^*\mu_2/p; r+1; p ) \\
\label{eq:c2cota2}
& = &\frac{ ( \Omega^*\mu_2/p-d-1  )^{(r)}}{(r
-1)!} \sum_{j=0}^\infty\frac{(-1)^j  ( \Omega^*\mu_2 /p-r-d-1
 )^{(j)}}{(r+j)j!} p^{r+j} \nonumber\\
% & = \frac{1}{(\nnum-1)!} \sum_{j=0}^\infty\frac{(-1)^j  (
& = & \frac{1}{(r-1)!} \sum_{j=0}^\infty\frac{(-1)^j (\Omega
^*\mu_2-(d+1)p) \cdots(\Omega^*\mu_2-(d+r+j)p)}{(r+j)j!}.
\end{eqnarray}
Similarly,
%
%e5.37 ###
\begin{equation}
\label{eq:c1cota2}
\tilde c_1(p) = \frac{1}{(r-1)!}
\sum_{j=0}^\infty\frac{(-1)^j (\Omega^*/\mu_1-dp) \cdots(\Omega
^*/\mu_1-(d+r+j-1)p)}{(r+j)j!}.
\end{equation}

According to (\ref{eq:c2cota2}) and (\ref{eq:c1cota2}), both
$\tilde c_1(p)$ and $\tilde c_2(p)$ can be expressed as power series
in~$p$: $\tilde c_1(p) = \sum_{i=0}^\infty u_i p^i$, $\tilde c_2(p) =
\sum_{i=0}^\infty v_i p^i$. Thus, the right-hand side of (\ref{eq:cc1c2opt}) is also a power series, $\sum_{i=0}^\infty w_i p^i$, with
$w_i = v_i - u_i$. The zero-order coefficient is
\[
w_0 = \frac{1}{(r-1)!}  \Biggl( \sum_{j=0}^\infty\frac{(-1)^j
(\Omega^*\mu_2)^{r+j}}{(r+j)j!} - \sum_{j=0}^\infty\frac
{(-1)^j (\Omega^*/\mu_1)^{r+j}}{(r+j)j!}  \Biggr).
\]
From the Taylor expansion of $\gamma(r,t)$ (\citet{Abramowitz70}, equation (6.5.29)),
\[
\gamma(r,t) = \frac{1}{\Gamma(r)} \sum_{j=0}^\infty\frac
{(-1)^j t^{r+j} }{(r+j)j!},
\]
the coefficient $w_0$ is recognized to be $c^*$, reflecting the fact
that the difference between the two sides of (\ref{eq:cc1c2opt})
is vanishingly small as $p \rightarrow0$. In addition, the first order
coefficient $w_1$ coincides with the derivative of the right-hand side
of (\ref{eq:cc1c2opt}) evaluated at $p=0$. It follows that \mbox{$w_1>0$}
is a~sufficient condition for the estimator (\ref{eq:hatpOmegaoptn+d}) to asymptotically guarantee the confidence level~$c^*$.%\looseness=1

The coefficient $w_1 = v_1 - u_1$ is computed as follows. The term
$v_1$ is obtained from (\ref{eq:c2cota2}) as%\looseness=1
%
%e5.38 ###
\begin{equation}
\label{eq:a21}
v_1 = \frac{1}{(r-1)!} \sum_{j=0}^\infty\frac
{ (-1)^{j+1} (\Omega^*\mu_2)^{r+j-1} [ (d+1)+\cdots+(d+r
+j) ] }
{ (r+j)j! }.
\end{equation}
Substituting $(d+1)+\cdots+(d+r+j) = (r+j) (r
+j+1+2d)/2$ into (\ref{eq:a21}) and taking into account that $\exp
(t) = \sum_{j=0}^\infty t^j/j!$, we have
\begin{eqnarray}
\label{eq:a21segunda}
v_1 & =& \frac{r+1+2d}{2(r-1)!} \sum_{j=0}^\infty\frac
{(-1)^{j+1} (\Omega^*\mu_2)^{r+j-1}}{j!} \nonumber
+ \frac{1}{2(r-1)!} \sum_{j=0}^\infty\frac{(-1)^j (\Omega \nonumber
^*\mu_2)^{r+j}}{j!}\qquad\\ \nonumber\\[-8pt]\\[-8pt]
& = & \frac{(\Omega^*\mu_2)^{r-1} \exp(-\Omega^*\mu_2 \nonumber
)}{2(r-1)!} (\Omega^*\mu_2-r-1-2d).\nonumber
\end{eqnarray}
Similarly,
%
%e5.39 ###
\begin{equation}
\label{eq:a11}
u_1 = \frac{(\Omega^*/\mu_1)^{r-1} \exp(-\Omega^*/\mu_1
)}{2(r-1)!} (\Omega^*/\mu_1-r+1-2d).
\end{equation}
From (\ref{eq:a21segunda}) and (\ref{eq:a11}), and making use
of (\ref{eq:Omega*funofdos}),
\begin{eqnarray*}
w_1 & = & \frac{(\Omega^*/\mu_1)^{r-1}}{2(r-1)!} [
M^{r-1} \exp(-\Omega^*\mu_2) (\Omega^*\mu_2-r-1-2d)
\\
&&\hphantom{\frac{(\Omega^*/\mu_1)^{r-1}}{2(r-1)!} [}
{}- \exp(-\Omega^*/\mu_1) (\Omega^*/\mu_1-r+1-2d)
] \\
% & = \frac{1}{2(\nnum-1)!}  (\frac{\nnum\log\fud}{\fud-1}
% )^{\nnum-1}  ( \fud^{\nnum-1}
% & \quad - \fud^{-\frac{\nnum}{\fud-1}}  (\frac{\nnum\log
% & = \frac{1}{2(\nnum-1)!}  (\frac{\nnum\log\fud}{\fud-1}  )^
% \fud^{-\frac{\nnum}{\fud-1}} \\
% & \quad\cdot ( \frac{1}{\fud}  (\frac{\nnum\fud\log\fud}{
% - \nnum-1-2d  ) - \frac{\nnum\log\fud}{\fud-1} + \nnum-1+2d
% )
& = & \frac{1}{2(r-1)!}  \biggl(\frac{r\log M}{M-1}
 \biggr)^{r-1}
M^{-r/(M-1)}  [ -(r+1+2d)/M+ r-1+2d ],
\end{eqnarray*}
which is positive if (\ref{eq:dmayorque}) holds. Thus, (\ref{eq:hatpOmegaoptn+d}) asymptotically guarantees the confidence $c^*$ for
$d$ as in (\ref{eq:dmayorque}).
\end{pf*}

\begin{lemma}\label{lemma:intsum}
For $r$, $n_0 \in\mathbb N$ with $r\geq3$ and
$
n_0 \leq(r- \sqrt{r}) / p,
$
%
%e5.40 ###
\begin{equation}
\label{eq:desiglemmaintsum}
\int_{r}^{n_0} t^{r-1} \exp(-pt) \,\mathrm{d}t \geq\sum
_{n=r}^{n_0-1} (n-1)^{(r-1)} (1-p)^{n
-r}.
\end{equation}
\end{lemma}

\begin{pf}
For $n_0$ as given, the sub-integral function in (\ref{eq:desiglemmaintsum}) is increasing within the integration range.
Consequently, for (\ref{eq:desiglemmaintsum}) to hold, it is
sufficient that
$
\label{eq:cond1}
n^{r-1} \exp(-np) \geq(n-1)^{(r-1)}
(1-p)^{n-r}
$
% for $\nden= \nnum, \ldots, \nden_0$, or, taking into account
for $n= r, \ldots, \lfloor(r- \sqrt{r})/p\rfloor
-1$. This condition can be shown to be satisfied by means of reasoning
analogous to
% vale
that in the proof of Lemma 1 of \citet{Mendo08a}, part~(i).\quad\mbox{}
\end{pf}

\begin{pf*}{Proof of Proposition \ref{prop:confgargen}}
The confidence for the estimator (\ref{eq:hatpOmegan+1}) is
expressed as $c_2(p) - c_1(p)$, $c_1(p) = P[N\leq n_1
-1]$, $c_2(p) = P[N\leq n_2]$, where $n_1$ and
$n_2$ are given by (\ref{eq:ndunonddos}) with $d=1$. Let $c'_2 =
1-c_2$. A similar argument as in the proof of Proposition \ref{prop:limdeconf}, based on the Poisson theorem, shows that $\lim_{p
\rightarrow0} c_1(p) = \bar c_1$ and $\lim_{p \rightarrow0} c'_2(p)
= \bar c'_2$, with $\bar c_1 = \gamma(r,\Omega/\mu_1)$, $\bar
c'_2 = 1-\gamma(r,\Omega\mu_2)$. Thus, to establish the desired
result, it suffices to prove that $c_1(p) < \bar c_1$ and $c'_2(p) <
\bar c'_2$ for $\mu_1$, $\mu_2$ as in (\ref{eq:condfunofdosgen}).

Regarding $c'_2(p)$, it is shown in Appendix C of \citet{Mendo06} that
$P[N\leq\lfloor a/p\rfloor] > \lim_{p\rightarrow0} P
[N\leq\lfloor a/p\rfloor]$ for $a > r+\sqrt r$.
Equivalently, for any $\theta> 1$,
\[
P\bigl[ N\geq\bigl\lfloor\theta\bigl(r+\sqrt r\bigr)/p\bigr\rfloor+ 1 \bigr]
< \lim_{p\rightarrow0} P\bigl[ N\geq\bigl\lfloor\theta\bigl(r
+\sqrt r\bigr)/p\bigr\rfloor+ 1 \bigr].
\]
This can be made to correspond to $c'_2(p) < \bar c'_2$ by taking
$\lfloor\Omega\mu_2/p-1\rfloor= \lfloor\theta(r+\sqrt r
)/p\rfloor$, that is,
%
%e5.41 ###
\begin{equation}
\label{eq:fdostheta}
\mu_2= \frac{\theta(r+\sqrt r)+p}{\Omega}.
\end{equation}
Thus, the inequality $c'_2(p) < \bar c'_2$ holds if $\mu_2$ is given
by (\ref{eq:fdostheta}) for some $\theta>1$ or, equivalently, if
$\mu_2> (r+\sqrt r+p)/\Omega$. Therefore, it holds for
arbitrary $p$ if $\mu_2$ satisfies the second inequality in~(\ref{eq:condfunofdosgen}).

As for $c_1(p)$, from (\ref{eq:f}), it can be expressed as
% As for $c_1(p)$, a similar technique to that used in \cite{Mendo08a}
%can be employed to establish that $c_1(p) < \bar c_1$. From \eqref{eq:
%f}, $c_1(p)$ is expressed as
%
%e5.42 ###
\begin{equation}
\label{eq:c1pdes}
c_1(p) = \sum_{n=r}^{n_1-1} f(n) = \frac{p^r
}{(r-1)!} \sum_{n=r}^{n_1-1} (n-1)^{(r -1)} (1-p)^{n-r}.
\end{equation}
According to (\ref{eq:ndunonddos}) with $d=1$,
%
%e5.43 ###
\begin{equation}
\label{eq:ndunomenor}
n_1< \frac{\Omega}{p\mu_1}
\end{equation}
and thus
%
%e5.45 ###
%e5.44 ###
\begin{eqnarray}\label{eq:barc1cota}
\bar c_1 &=& \gamma( r, \Omega/\mu_1)
% = \frac{1}{(\nnum-1)!} \int_{0}^{ \Omega/\funo} t^{\nnum-1} \exp(-t)
 = \frac{p^r}{(r-1)!} \int_{0}^{ \Omega/(p\mu_1) }
t^{r-1} \exp(-pt) \,\mathrm{d}t \nonumber\\[-8pt]\\[-8pt]
& >& \frac{p^r}{(r-1)!} \int_{r}^{n_1} t^{r-1}
\exp(-pt) \,\mathrm{d}t.\nonumber
\end{eqnarray}
From (\ref{eq:c1pdes}), (\ref{eq:barc1cota}) and Lemma \ref
{lemma:intsum}, it follows that the inequality $c_1(p) < \bar c_1$ is
satisfied if
%
%e5.46 ###
\begin{equation}
\label{eq:nduno1leq}
n_1\leq\frac{r-\sqrt r}{p}.
\end{equation}
Using (\ref{eq:ndunomenor}), it is seen that (\ref{eq:nduno1leq}) is fulfilled if $\mu_1$ satisfies the first inequality in (\ref{eq:condfunofdosgen}).
\end{pf*}

\begin{pf*}{Proof of Theorem \ref{theo:confgaropt}}
For $\Omega=\Omega^*$, using (\ref{eq:fud}) and (\ref{eq:Omega*funofdos}), the conditions in (\ref{eq:condfunofdosgen}) are
written as (\ref{eq:condoptfud34}) and (\ref{eq:condoptfud5omas}). The left-hand sides of (\ref{eq:condoptfud34}) and (\ref
{eq:condoptfud5omas}) are increasing functions of $M$, whereas
the right-hand sides decrease with $r$. These inequalities can
thus be written as $M\geq h_1(r)$, $M\geq h_2(r)$,
where $h_1$, $h_2$ are decreasing functions. This proves that for each
$r$, one of the inequalities implies the other. Furthermore,
defining $h(r)=\min\{h_1(r),h_2(r)\}$, the allowed
range for $M$ is expressed as $M\geq h(r)$ and $h$ is a
decreasing function. It only remains to prove that the limiting
condition is (\ref{eq:condoptfud34}) for $r=\{3, 4\}$ and
(\ref{eq:condoptfud5omas}) for $r\geq5$.

The left-hand sides of (\ref{eq:condoptfud34}) and (\ref{eq:condoptfud5omas}) are continuous functions of $M>1$.
Considering $r$ as if it were a continuous variable, the
right-hand sides are also continuous functions.
% vale "also" ah para indicar que se aplica al sujeto.
Assume that (\ref{eq:condoptfud34}) implies (\ref{eq:condoptfud5omas}) for a given $r_1$, that is, that the latter is
satisfied when the former holds with equality. Likewise, assume that
(\ref{eq:condoptfud5omas}) implies (\ref{eq:condoptfud34})
for a given $r_2$. The continuity of the involved functions then
implies that there exists $t \in[\min\{r_1,r_2\},\max\{
r_1,r_2\}]$ such that both (\ref{eq:condoptfud34}) and
(\ref{eq:condoptfud5omas}) hold with equality for $r=t$,
that is,
%
%e5.47 ###
\begin{equation}
\label{eq:condoptfudx}
\frac{M-1}{\log M} = \frac{t + \sqrt t}{t-1}, \qquad
\frac{M\log M}{M-1} = \frac{t + \sqrt t + 1}{t}.
\end{equation}
Multiplying both equalities in (\ref{eq:condoptfudx}) and
substituting into the first yields
%
%e5.48 ###
\begin{equation}
\label{eq:condoptxigual}
\log\frac{t+\sqrt t + 1}{t-\sqrt t} - \frac{2\sqrt t +1}{t} = 0.
\end{equation}
It is easily shown that (\ref{eq:condoptxigual}) has only one
solution, which lies in the interval $(4,5)$.
% Thus the limiting condition is the same for $\nnum_1$ and $\nnum_2$
%if $\nnum_1, \nnum_2 \in\{3,4\}$ or if $\nnum_1, \nnum_2 \in\{5,6,7,
Thus, one of the two conditions (\ref{eq:condoptfud34}) and
(\ref{eq:condoptfud5omas}) is the limiting one for $r\in\{
3,4\}$, whereas the other is for $r=\{5,6,7,\ldots\}$.
Taking any value from each set, the limiting condition is seen to be
(\ref{eq:condoptfud34}) in the former case and (\ref{eq:condoptfud5omas}) in the latter.
\end{pf*}

\section*{Acknowledgements}
The authors wish to thank The Editor, Prof. Holger Rootz\'en, and one referee for their positive
comments. The first author also wishes to thank Dr. Manuel Kauers and
Prof. Wolfram Koepf for help with certain expressions involving
Stirling numbers and hypergeometric terms, respectively, and Dr. Xinjia
Chen for reading an early version of this paper.

\printhistory

\end{document}